\documentclass[11pt]{amsart} 
\usepackage[curve,v2]{xy}

\usepackage{amssymb,latexsym,amsfonts}

\newcommand{\tensor}{\otimes}

\newcommand{\rk}{\operatorname{rank}}
\newcommand{\Hom}{\operatorname{Hom}}

\renewcommand{\c}{\operatorname{c}}

\newcommand{\picard}{\operatorname{Pic}}
\newcommand{\Ext}{\operatorname{Ext}}
\newcommand{\ext}{\operatorname{ext}}
\newcommand{\sext}{\operatorname{\mathcal{E}\mathit{xt}}}

\newcommand{\shom}{\operatorname{\mathcal{H}\mathit{om}}}
\renewcommand{\H}{\operatorname{H}}
\newcommand{\h}{\operatorname{h}}

\newcommand{\wt}[1]{\widetilde{#1}}

\newcommand{\blow}[2]{{\rm Bl}_{#2}{(#1})}

\newcommand{\ses}[3]{0\rightarrow#1\rightarrow#2
   \rightarrow#3\rightarrow0}

\newcommand{\M}{{\mathcal M}}

\newcommand{\F}{{\mathcal F}}
\newcommand{\G}{{\mathcal G}}

\newcommand{\I}{{\mathcal I}}
\renewcommand{\L}{{\mathcal L}}
\renewcommand{\O}{{\mathcal O}}
\renewcommand{\P}{{\mathbb{P}}}

\newcommand{\Z}{{\mathbb{Z}}}

\renewenvironment{proof}{\par \medskip \noindent
{\sc Proof:}}{}

\begingroup
\newtheorem{thm}{Theorem}   
\newtheorem{cor}[thm]{Corollary}     
\newtheorem{lemma}[thm]{Lemma}         
\newtheorem{prop}[thm]{Proposition}  

\newtheorem{defn}[thm]{Definition}   
   
\newtheorem{notation}[thm]{Notation}

\endgroup

\newenvironment{rem}[2]{\refstepcounter{thm} \label{#2} 
\par \medskip \noindent {\bf #1 \thethm .}}{\par \medskip}

\begin{document}

\pagenumbering{arabic}

\title[Moduli of Reflexive Sheaves]{Moduli of Reflexive Sheaves on Smooth Projective 3-folds}

\author[Peter Vermeire]{Peter Vermeire}

\address{Department of Mathematics, 214 Pearce, Central Michigan
University, Mount Pleasant MI 48859}

\email{verme1pj@cmich.edu}

\date{\today}

\begin{abstract} 
We compute the expected dimension of the moduli space of torsion-free rank 2 sheaves at a point corresponding to a stable reflexive sheaf and give conditions for the existence of a perfect tangent-obstruction complex on a class of smooth projective threefolds; this class includes Fano and Calabi-Yau threefolds.  We also explore both local and global relationships between moduli spaces of reflexive rank 2 sheaves and the Hilbert scheme of curves.
\end{abstract}

\maketitle

\section{Introduction and Preliminaries}
We work over an algebraically closed field of characteristic $0$. 

In this paper, motivated by Qin's \cite{qin92}, \cite{qin1} and Qin and Li's \cite{lq95} work on the relationship between the moduli of vector bundles on surfaces and Hilbert schemes of points, by Hartshorne's work
\cite{hartvb},\cite{hart},\cite{hart2},\cite{hart3} on curves in
$\P^3$ via the Serre Correspondence, and by the theory of virtual fundamental classes and Donaldson-Thomas invariants \cite{Thomas}, we study the moduli of reflexive rank $2$ sheaves on smooth projective threefolds.  Some very precise statements concerning the structure of these moduli spaces on Fano threefolds (especially smooth hypersurfaces) can be found in \cite{ball}, \cite{ball2}, \cite{beau}, \cite{druel}, \cite{ilman}, \cite{IM}, \cite{IM2}, \cite{IM3}, \cite{IR}, \cite{MT}, \cite{T}.

Recall that a  coherent sheaf $\F$ is \textit{torsion-free} if the natural map of $\F$ to its double-dual $h:\F\rightarrow\F^{**}$ is injective, and that $\F$ is \textit{reflexive} if $h$ is an isomorphism. 
We refer the reader to \cite{hart} for basic properties of reflexive
sheaves.  Recall the following {\it Serre Correspondence} for
reflexive sheaves:
\begin{thm}\label{scorr}\cite[4.1]{hart}
Let $X$ be a smooth projective threefold, $M$ an invertible sheaf
with $\H^1(X,M^*)=\H^2(X,M^*)=0$.  There is a one-to-one correspondence
between 
\begin{enumerate}
\item pairs $(\F,s)$ where $\F$ is a rank $2$ reflexive sheaf on $X$
  with $\det\F=M$ and $s\in\Gamma(\F)$ is a section whose 
  set has codimension $2$
\item pairs $(Y,\xi)$ where $Y$ is a closed Cohen-Macaulay curve in
  $X$, generically a local complete intersection, and
  $\xi\in\Gamma(Y,\omega_Y\tensor\omega_X^*\tensor M^*)$ is a section
  which generates the sheaf $\omega_Y\tensor\omega_X^*\tensor M^*$
  except at finitely many points.
\end{enumerate}
Furthermore, $\c_3(\F)=2p_a(Y)-2+\c_1(X)\c_2(\F)-\c_1(\F)\c_2(\F)$.
\nopagebreak \hfill $\Box$ \par \medskip
\end{thm}

Note that if $\F$ is locally free, then the corresponding curve $Y$ is a local complete
intersection,  $\omega_Y\tensor\omega_X^*\tensor M^*\cong\O_Y$, $\xi$
is a non-zero section, and $c_3(\F)=0$.  In this case we say $Y$ is
{\em subcanonical}.

Note that any smooth curve in $X$ can be made subcanonical by blowing up:
\begin{prop}
Let $X$ be a smooth projective variety of dimension $n$, $C\subset X$ a smooth curve.  Then there exists a finite set $D\subset C$ and a line bundle $L\in\picard(X)$ such that if $\pi:\wt{X}=\blow{X}{D}\rightarrow X$ then
\begin{enumerate}
\item $\omega_{\wt{C}}\tensor \pi^*L^*\tensor\omega_{\wt{X}}^*=\O_{\wt{C}}$, and so the proper transform of $C$ is subcanonical.
\item $\H^i(\wt{X},\pi^*L^*)=0$ for $i\geq1$.
\end{enumerate}
In particular, if $n=3$ there is a rank two locally free sheaf $\F$ on $\wt{X}$ with $\wedge^2\F=\pi^*L$ and with a section $s\in\H^0(\wt{X},\F)$ whose zero scheme is $\wt{C}$. 
\end{prop}

\begin{proof}
Let $L\in\picard(X)$ be such that 
\begin{enumerate}
\item $\H^i(X,L^*)=0$ for $i\geq1$.
\item $\omega_C\tensor L^*\tensor\omega_X^*=\O_C((n-1)D)$ where $D\subset C$ is a collection of distinct points.
\end{enumerate}
Both conditions are satisfied by a sufficiently ample $L^*$.

Now, $\H^i(\wt{X},\pi^*L^*)=0$ for $i\geq1$ follows immediately.  Letting $\{E_i| i=1,\ldots,|D|\}$ be the exceptional divisors we have 
\begin{eqnarray*}
\omega_{\wt{C}}\tensor \pi^*L^*\tensor\omega_{\wt{X}}^*&=&\omega_{\wt{C}}\tensor \pi^*\left(L^*\tensor\omega_X^*\right)\tensor\O_{\wt{C}}((1-n)E_i)\\
&=&\omega_C\tensor L^*\tensor\omega_X^*\tensor\O_C((1-n)D)\\
&=&\O_C 
\end{eqnarray*}
\nopagebreak \hfill $\Box$ \par \medskip
\end{proof}

Unfortunately from the point of view of moduli, in the case of a threefold the associated locally free sheaf $\F$ on $\wt{X}$ is in some sense less likely to be stable than the original reflexive sheaf on $X$ (though perhaps this could be repaired by considering the moduli of pairs).  Thus we do not pursue this direction and instead in Section~\ref{maincalc} we study general properties of the moduli spaces of reflexive sheaves on $X$.  In particular, we calculate (Corollary~\ref{myfano}) the expected dimension of the space in cases where the canonical bundle is at least somewhat non-positive.  Further, we show that the space of stable, rank $2$ reflexive sheaves admits a perfect tangent-obstruction complex (Corollary~\ref{perfect}) in a significant class of varieties not covered by Thomas' results \cite{Thomas}.  

In Section~\ref{localtheory}, we study of the local relationship between the moduli space of reflexive sheaves, the moduli space of ideal sheaves, and the Hilbert scheme of curves (Theorem~\ref{bigthm}).  Section~\ref{globaltheory} is concerned with the global relationship between these moduli spaces; we give some results (Proposition~\ref{hilbdominates}) in the case of vector bundles of rank $2$.

We recall some basic terms and results:

\begin{defn}
Let $L$ be an ample line bundle on a smooth projective variety $X$. A torsion-free sheaf $\F$ on $X$ is {\bf \boldmath $L$-stable} (resp. {\bf \boldmath $L$-semistable}) if for every coherent subsheaf $\F'$ of $\F$ with $0<\rk\F'<\rk\F$, we have $\mu(\F',L)<\mu(\F,L)$ (resp. $\leq$), where for any coherent sheaf $\G$ we define
$$\mu(\G,L)=\frac{\c_1(\G).[L]^{\dim X-1}}{\left( \rk\G\right) [L]^{\dim X}}$$
Note that if $\rk\F=2$, it suffices to take $\F'$ invertible.  If $X\subset\P^n$ we simply say $\F$ is (semi)stable to mean $\F$ is $\O_X(1)$-(semi)stable.  A coherent sheaf $\F$ on a projective variety $X$ is called {\bf simple} if $\dim\Hom_{\O_X}(\F,\F)=1$.  Recall (\cite[1.2.8]{lehnbook}) that stable sheaves are simple.
\nopagebreak \hfill $\Box$ \par \medskip
\end{defn}

\begin{thm}\label{rr}\cite[A.5.3]{hartshorne}
Let $\F$ be a coherent sheaf of rank $r$ on a smooth projective threefold $X$.  The Riemann-Roch formula is
\begin{eqnarray*}
\chi(X,\F) & = &
\frac{1}{6}\c_1^3(\F)-\frac{1}{2}\c_1(\F)\c_2(\F)-\frac{1}{2}\c_1(X)\c_2(\F)+\frac{1}{4}\c_1(X)\c_1^2(\F)
\\
 & & +\frac{1}{12}\c_1^2(X)\c_1(\F)+\frac{1}{12}\c_2(X)\c_1(\F)+\frac{r}{24}\c_1(X)\c_2(X)+\frac{1}{2}\c_3(\F)
\end{eqnarray*}
Note that if $\F$ is a locally free sheaf of rank $r$, then
$$\chi(\F)+\chi(\F^*)=-\c_1(X)\c_2(\F)+\frac{1}{2}\c_1(X)\c_1^2(\F)+\frac{r}{12}\c_1(X)\c_2(X)$$
\nopagebreak \hfill $\Box$ \par \medskip
\end{thm}

\begin{prop}\label{localtoglobal}\cite[2.5]{hart}
Let $\F$ be a reflexive sheaf on a normal projective threefold $X$, $\G$ a sheaf of $\O_X$-modules.  Then there are isomorphisms

$\H^0(X,\sext^0_{\O_X}(\F,\G))=\Ext^0_{\O_X}(\F,\G)$

$\H^3(X,\sext^0_{\O_X}(\F,\G))=\Ext^3_{\O_X}(\F,\G)$

and an exact sequence
$$0\rightarrow \H^1(X,\sext^0_{\O_X}(\F,\G))\rightarrow \Ext^1_{\O_X}(\F,\G)\rightarrow \H^0(X,\sext^1_{\O_X}(\F,\G))$$
$$\rightarrow \H^2(X,\sext^0_{\O_X}(\F,\G))\rightarrow \Ext^2_{\O_X}(\F,\G)\rightarrow0$$
\nopagebreak \hfill $\Box$ \par \medskip
\end{prop}

\section{Extension Calculations}\label{maincalc}

\begin{notation}
For a coherent sheaf $\F$ of rank $r$, we denote $\Delta(\F)=2r\c_2(\F)-(r-1)\c_1^2(\F)=\c_2(\F\tensor\F^*)$.  Note that $\Delta(\F)=\Delta(\F\tensor L)$ for any invertible sheaf $L$.
\nopagebreak \hfill $\Box$ \par \medskip
\end{notation}

It is known \cite{mar},\cite{mar2} that there is a coarse projective moduli space for semistable torsion-free sheaves with given Chern classes on a smooth projective threefold $X$.  The tangent space at a stable sheaf $\F$ is $\Ext^1_{\O_X}(\F,\F)$, the obstructions lie in $\Ext^2_{\O_X}(\F,\F)$, and if $\Ext^2_{\O_X}(\F,\F)=0$ then the moduli space is smooth at $\F$.
This motivates:
\begin{defn}
Let $\F$ be a torsion-free sheaf on a smooth projective threefold $X$.  The \textbf{expected dimension} of the coarse moduli space of stable torsion-free sheaves with Chern classes equal to that of $\F$ is $$\mathfrak{D}(\F)=\dim\Ext^1_{\O_X}(\F,\F)-\dim\Ext^2_{\O_X}(\F,\F)$$
\nopagebreak \hfill $\Box$ \par \medskip
\end{defn} 

We make a formal computation:
\begin{prop}\label{bigextcomp}
Let $\F$ be a rank $r$ coherent sheaf of homological dimension $1$ on a smooth projective threefold $X$.  Then $$ \sum_{i=0}^3(-1)^i\dim\Ext^i_{\O_X}(\F,\F)=\frac{r^2\c_1(X)\c_2(X)}{24}-\frac{\c_1(X)}{2}\Delta(\F)$$
\end{prop}

\begin{proof}
By definition, because $\F$ has homological dimension $1$ there is a resolution 
$$\ses{E_1}{E_0}{\F}$$
with $E_0$, $E_1$ locally free of rank $k+r$ and $k$ respectively.  Furthermore, as in \cite[3.4]{hart} we compute the global extension groups as the hypercohomology of the complex
$$E_0^*\tensor E_1\rightarrow \left(E_0^*\tensor E_0\right)\oplus \left(E_1^*\tensor E_1\right)\rightarrow E_1^*\tensor E_0$$ and so
 $$\sum_{i=0}^3(-1)^i\dim\Ext^i(\F,\F)=\chi(E_0^*\tensor E_0)+\chi(E_1^*\tensor E_1)-\chi(E_0^*\tensor E_1)-\chi(E_1^*\tensor E_0)$$
In general, if $r_1=\rk(E)$ and $r_2=\rk(F)$, one finds
\begin{eqnarray*}
\c_1(E\tensor F)&=&r_2\c_1(E)+r_1\c_1(F)\\
\c_2(E\tensor F) &=&\binom{r_2}{2}\c_1^2(E)+r_2\c_2(E)+ (r_1r_2-1)\c_1(E)\c_1(F)\\
&&+r_1\c_2(F)+\binom{r_1}{2}\c_1^2(F)\\
c_3(E\tensor F) &=&\binom{r_2}{3}\c_1^3(E)+2\binom{r_2}{2}\c_1(E)\c_2(E)+(r_1r_2-2)\c_1(E)\c_2(F)\\
&&+\frac{1}{2}(r_2-1)(r_1r_2-2)\c_1^2(E)\c_1(F)+\frac{1}{2}(r_1-1)(r_1r_2-2)\c_1(E)\c_1^2(F)\\
&&+(r_1r_2-2)\c_2(E)\c_1(F)+2\binom{r_1}{2}\c_1(F)\c_2(F)+\binom{r_1}{3}\c_1^3(F)\\
&&+r_2(r_2^2-3r_2+3)\c_3(E)+r_1(r_1^2-3r_1+3)\c_3(F)
\end{eqnarray*}
Substituting in our case, we obtain: 
\begin{eqnarray*}
\c_1(E_i\tensor E_i^*)&=&\c_3(E_i\tensor E_i^*)=0\\
\c_2(E_0\tensor E_0^*)&=&2(k+r)\c_2(E_0)-(k+r-1)\c_1^2(E_0)=\Delta(E_0)\\
\c_2(E_1\tensor E_1^*)&=&2k\c_2(E_1)-(k-1)\c_1^2(E_1)=\Delta(E_1)
\end{eqnarray*}

Computing via Riemann-Roch (Theorem~\ref{rr}) we see
$$\chi(E_0\tensor E_0^*)=-\frac{1}{2}\c_1(X)\left(2(k+r)\c_2(E_0)-(k+r-1)\c_1^2(E_0)-\frac{(k+r)^2}{12}\c_2(X)\right)$$
$$\chi(E_1\tensor E_1^*)=-\frac{1}{2}\c_1(X)\left(2k\c_2(E_1)-(k-1)\c_1^2(E_1)-\frac{k^2}{12}\c_2(X)\right)$$

For the other two terms, we again substitute into the general formulae to obtain:
\begin{eqnarray*}
c_1(E_0^*\tensor E_1)&=&(k+r)c_1(E_1)-kc_1(E_0)\\
c_2(E_0^*\tensor E_1)&=&\binom{k}{2}c_1^2(E_0)+kc_2(E_0)-(k^2+rk-1)c_1(E_0)c_1(E_1)\\
&&+(k+r)c_2(E_1)+\binom{k+r}{2}c_1^2(E_1)
\end{eqnarray*}

and so by Theorem~\ref{rr} we have
$$\chi(E_0^*\tensor E_1)+\chi(E_0\tensor E_1^*)=\frac{c_1(X)}{2}\left(c_1^2(E_0^*\tensor E_1)+\frac{k(k+r)}{6}c_2(X)-2c_2(E_0^*\tensor E_1)\right)$$

Putting the four terms together:
\begin{eqnarray*}
\sum_{i=0}^3(-1)^i\dim\Ext^i(\F,\F)&=&\chi(E_0^*\tensor E_0)+\chi(E_1^*\tensor E_1)-\chi(E_0^*\tensor E_1)-\chi(E_1^*\tensor E_0)\\
&=&\frac{r^2c_1(X)c_2(X)}{24}+rc_1(X)\left(c_2(E_1)-c_2(E_0)\right)\\
&&+\frac{c_1(X)}{2}\left((r-1)c_1^2(E_0)-(r+1)c_1^2(E_1)+2c_1(E_0)c_1(E_1)\right)\\
&=&\frac{r^2c_1(X)c_2(X)}{24}-\frac{c_1(X)}{2}\Delta(\F)
\end{eqnarray*}
Where the last equality holds after replacing $c_1(\F),c_2(\F)$ with the appropriate Chern classes of $E_0$ and $E_1$ as derived from the resolution of $\F$.
\nopagebreak \hfill $\Box$ \par \medskip
\end{proof}

\begin{cor}
Let $\F$ be a rank $2$ reflexive sheaf on a smooth projective threefold $X$.  Then $$ \sum_{i=0}^3(-1)^i\dim\Ext^i_{\O_X}(\F,\F)=\frac{c_1(X)c_2(X)}{6}-\frac{c_1(X)}{2}\Delta(\F)$$
\end{cor}

\begin{proof}
This follows immediately from Proposition~\ref{bigextcomp} noting that reflexive sheaves have homological dimension $1$ \cite[1.2]{hart}.
\nopagebreak \hfill $\Box$ \par \medskip
\end{proof}

Since stable sheaves are simple, to compute $\mathfrak{D}(\F)$ for a stable sheaf we need only compute $\dim\Ext^3_{\O_X}(\F,\F)$.  We give two results, one for varieties with effective anticanonical divisor (Proposition~\ref{killext3fano}) and one for a class including Fano varieties (Proposition~\ref{weakerthanfano}).

\begin{prop}\label{killext3fano}
Let $\F$ be a reflexive sheaf on a smooth projective threefold $X$ with $\omega_X^*$ effective.  
\begin{enumerate}
\item If $\omega_X=\O_X$ then $\Ext^3_{\O_X}(\F,\F)^*=\Hom_{\O_X}(\F,\F)$.
\item If $\omega_X\neq\O_X$ then $\dim\Ext^3_{\O_X}(\F,\F)<\dim\Hom_{\O_X}(\F,\F)$.
\end{enumerate}
\end{prop}

\begin{proof}(Cf. \cite[3.39]{Thomas})
We have 
\begin{eqnarray*}
\Ext^3_{\O_X}(\F,\F)&=&\H^3(X,\shom_{\O_X}(\F,\F))\\
&=&\Hom_{\O_X}(\shom_{\O_X}(\F,\F),\omega_X)^*
\end{eqnarray*} 
where the first equality is Proposition~\ref{localtoglobal} and the second is Serre Duality.  We also have
\begin{eqnarray*}
\Hom_{\O_X}(\F,\F)&=&\H^0(X,\shom_{\O_X}(\F,\F))\\
&=&\H^3(X,\omega_X\tensor\shom_{\O_X}(\F,\F))^*\\
&=&\Hom_{\O_X}(\omega_X\tensor\shom_{\O_X}(\F,\F),\omega_X)\\
&=&\Hom_{\O_X}(\shom_{\O_X}(\F,\F),\O_X)
\end{eqnarray*} 
where the first equality is Proposition~\ref{localtoglobal}; the second is \cite[2.5]{hart}; the third is Serre Duality; the fourth is \cite[III.6.7]{hartshorne}.

Now clearly if $\omega_X=\O_X$ then $\Ext^3_{\O_X}(\F,\F)^*=\Hom_{\O_X}(\F,\F)$.  Otherwise, letting $\omega_X^*=\O_X(D)$ we have 
$$\ses{\omega_X}{\O_X}{\O_D}$$
Applying $\Hom_{\O_X}(\shom_{\O_X}(\F,\F),\cdot)$ yields
$$0\rightarrow \Ext^3_{\O_X}(\F,\F)^*\rightarrow
\Hom_{\O_X}(\F,\F)\xrightarrow{f}
\Hom_{\O_X}(\shom_{\O_X}(\F,\F),\O_D)\rightarrow\cdots$$
where $f$ is not the zero map as it preserves at least the homotheties of $\F$.
\nopagebreak \hfill $\Box$ \par \medskip
\end{proof}

\begin{cor}
Let $\F$ be a reflexive sheaf on a smooth projective threefold $X$ with $\omega_X=\O_X$.  Then $\mathfrak{D}(\F)=0$
\end{cor}

\begin{proof}
As $\c_1(X)=\c_1(\O_X)$, this follows immediately from Propositions~\ref{bigextcomp} and ~\ref{killext3fano}.
\nopagebreak \hfill $\Box$ \par \medskip
\end{proof}

We also have the following, which applies, in particular, to Fano varieties:
\begin{prop}\label{weakerthanfano}
Let $\F$ be a stable rank two reflexive sheaf on a smooth projective threefold $X\subset\P^n$.  If there exists an $n\in\Z$ such that $\H^0(X,\F\tensor\omega_X^n)\neq0$ and $\H^0(X,\F\tensor\omega_X^{n+1})=0$, then $\Ext^3_{\O_X}(\F,\F)=0$.
\end{prop}

\begin{proof}
By hypothesis, for some $n$ we have $\H^0(X,\F\tensor\omega_X^n)\neq0$ and $\H^0(X,\F\tensor\omega_X^{n+1})=0$.  Let $D\subset X$ be an effective divisor such that $\F\tensor\omega_X^n\tensor\O_X(-D)$ has a section whose zero locus is a curve $C$; note that $D$ may be empty and that $\H^0(X,\F\tensor\omega_X^{n+1}\tensor\O_X(-D))=0$.  Letting $\G=\F\tensor\omega_X^n\tensor\O_X(-D)$, we compute $\Ext^3_{\O_X}(\G,\G)$.

The section of $\G$ induces an exact sequence
$$\ses{\O_X}{\G}{\I_C\tensor\det\G}$$

Applying $\Hom_{\O_X}(\cdot,\G)$ yields
$$\cdots\rightarrow\Ext^3_{\O_X}(\I_C\tensor\det\G,\G)\rightarrow\Ext^3_{\O_X}(\G,\G)\rightarrow \H^3(X,\G)\rightarrow0$$

We have $\H^3(X,\G)=\H^0(X,\G^*\tensor\omega_X)^*= \H^0(X,\G\tensor\omega_X\tensor\det\G^*)^*$ but the last group is zero by stability of $\G$ and the above discussion. 

Applying $\Hom_{\O_X}(\cdot,\G)$ to the basic sequence $$\ses{\I_C\tensor\det\G}{\det\G}{\O_C\tensor\det\G}$$
yields $$\cdots\rightarrow\Ext^3_{\O_X}(\det\G,\G)\rightarrow\Ext^3_{\O_X}(\I_C\tensor\det\G,\G)\rightarrow 0$$
We know $\Ext^3_{\O_X}(\det\G,\G)=\H^3(X,\G^*)=\H^0(X,\G\tensor\omega_X)=0$, therefore $\Ext^3_{\O_X}(\I_C\tensor\det\G,\G)=0$ and the vanishing of $\Ext^3_{\O_X}(\G,\G)$ follows.
\nopagebreak \hfill $\Box$ \par \medskip
\end{proof}

\begin{cor}\label{myfano}
Let $\F$ be a stable rank two reflexive sheaf on a smooth projective threefold $X$.  Assume either that $\omega_X^*$ is non-trivial and effective or that there exists an $n\in\Z$ such that $\H^0(X,\F\tensor\omega_X^n)\neq0$ and $\H^0(X,\F\tensor\omega_X^{n+1})=0$.  Then
$$\mathfrak{D}(\F)=1-\frac{\c_1(X)\c_2(X)}{6}+\frac{\c_1(X)\Delta(\F)}{2}$$
\nopagebreak \hfill $\Box$ \par \medskip
\end{cor}

It is interesting to see when the expected dimension is zero in the easily computable cases of smooth Fano complete intersection threefolds.  In particular, if $X$ is a hypersurface of degree $r\leq 4$ in $\P^4$ then a case-by-case examination shows that $\mathfrak{D}(\F)$ cannot be $0$ if $r=1,3,4$.  Similarly, if $X$ is a smooth complete intersection $X\subseteq\P^5$ of type $(2,2)$ or is a smooth complete intersection $X\subseteq\P^6$ of type $(2,2,2)$, then $\mathfrak{D}(\F)$ cannot be $0$.  This leaves two cases which serve as useful examples in Section~\ref{globaltheory}.

\begin{rem}{Example}{luckwithr=2}
Let $\F$ be a stable rank two reflexive sheaf on a smooth quadric hypersurface $X\subseteq\P^4$ and suppose that $\wedge^2\F=\O_X(k)$.  Then $\mathfrak{D}(\F)=0$ if and only if $k^2+1=2\c_1(\O_X(1))\c_2(\F)$.

Let $C\subset X$ be a line.  We associate to $C$ a rank $2$ vector bundle $\F$ with $\wedge^2\F=\O_X(1)$.  One can show that $\F$ is stable and that $\Ext^2_{\O_X}(\F,\F)=0$, hence the moduli space is smooth of dimension $0$.  Note, however, that $\h^0(C,N_{C/X})=3$ and $\H^1(C,N_{C/X})=0$, hence the Hilbert scheme is smooth of dimension three at $C$.
\nopagebreak \hfill $\Box$ \par \medskip
\end{rem}

\begin{rem}{Example}{luckwith(2,3)}
Let $\F$ be a stable rank two reflexive sheaf on a smooth complete intersection $X\subseteq\P^5$ of type $(2,3)$.   Then $\c_1(X)=\c_1(\O_X(1))$ and $\c_2(X)=4\c_1^2(\O_X(1))$.  If $\wedge^2\F=\O_X(k)$ then $\mathfrak{D}(\F)=0$ exactly when $2\c_1(\O_X(1))\c_2(\F)=3(1+k^2)$.

Let $C$ be a smooth plane cubic and $X$ a smooth complete intersection of type $(2,3)$ which contains it.
We associate to $C$ a rank $2$ vector bundle $\F$ with $\wedge^2\F=\O_X(1)$.  One can show that $\F$ is stable and that $\Ext^2_{\O_X}(\F,\F)=0$, hence the moduli space is again smooth of dimension $0$.  Note, however, that $\h^0(C,N_{C/X})=3$ and $\H^1(C,N_{C/X})=0$, hence the Hilbert scheme is again smooth of dimension three at $C$.
\nopagebreak \hfill $\Box$ \par \medskip
\end{rem}

By \cite[3.30]{Thomas} and \cite[3.7]{LiTian}, Proposition~\ref{weakerthanfano} also immediately implies:

\begin{cor}\label{perfect}
Let $X$ be a smooth projective threefold and let $\M$ be the moduli space of rank $2$ semistable sheaves with Chern classes $\c_i$ and determinant $L$.  Suppose that all such sheaves are stable and reflexive and suppose that for each $\F\in\M$ there exists an $n\in\Z$ such that $\H^0(X,\F\tensor\omega_X^n)\neq0$ and $\H^0(X,\F\tensor\omega_X^{n+1})=0$ (e.g. if $X$ is Fano).  Then $\M$ admits a perfect tangent-obstruction complex; further, there is a virtual cycle $Z_0\subset\M$ of dimension $\mathfrak{D}(\F)$ defined by the tangent-obstruction functors.
\nopagebreak \hfill $\Box$ \par \medskip
\end{cor}

\section{The Local Structure of Moduli of Reflexive Sheaves}\label{localtheory}

Several results in this section contain the hypothesis that $\H^2(X,\F)=0$.  To show this is not very restrictive (especially when $\F$ is stable and $X$ is Fano) we have:

\begin{lemma}\label{F}
Let $X$ be a smooth projective threefold with $\H^2(X,\O_X)=0$, $\F$ a rank $2$ reflexive sheaf with $\H^2(X,\det\F)=0$.  Suppose that $\F$ has a section whose zero scheme is a curve $C$ and that $\det\F\tensor\O_C$ is non-special.  Then $\H^2(X,\F)=0$.
\end{lemma}

\begin{proof}
This follows immediately from the sequence $\ses{\O_X}{\F}{\I_C\tensor\det\F}$.
\nopagebreak \hfill $\Box$ \par \medskip
\end{proof}

We compare the local structure of moduli of reflexive sheaves to moduli of ideal sheaves (see \cite{Vakil} and the remarks after Proposition~\ref{hilbdominates} for relations between this and the Hilbert scheme) 

\begin{prop}\label{betterhilbvsreflex}
Let $X$ be a smooth projective threefold, $\F$ a rank $2$ reflexive sheaf.
Suppose that $\F$ has a section whose zero scheme is a curve $C$.
\begin{enumerate}
\item If $\H^2(X,\F)=0$ and $\H^1(X,\I_C\tensor\det\F\tensor\omega_X)=0$, then the vanishing of $\Ext^2_{\O_X}(\I_C,\I_C)$ implies the vanishing of $\Ext^2_{\O_X}(\F,\F)$.
\item If $\H^0(X,\F\tensor\omega_X)=0$ and $\H^1(X,\I_C\tensor\det\F)=0$, then the vanishing of $\Ext^2_{\O_X}(\F,\F)$ implies the vanishing of $\Ext^2_{\O_X}(\I_C,\I_C)$.
\end{enumerate}
\end{prop}

\begin{proof}

\underline{Part (1)}:

Applying $\Hom_{\O_X}(\cdot,\F)$ to the sequence
$$\ses{\O_X}{\F}{\I_C\tensor\det\F}$$
we have 
$$\cdots\rightarrow\Ext^2_{\O_X}(\I_C\tensor\det\F,\F)\rightarrow \Ext^2_{\O_X}(\F,\F)\rightarrow \Ext^2_{\O_X}(\O_X,\F)\rightarrow\cdots$$
where $\Ext^2_{\O_X}(\O_X,\F)=\H^2(X,\F)=0$ by hypothesis.  Therefore, it suffices to show that $\Ext^2_{\O_X}(\I_C\tensor\det\F,\F)=\Ext^2_{\O_X}(\I_C,\F^*)=0$.

Applying $\Hom_{\O_X}(\I_C,\cdot)$ to the sequence
$$\ses{\det\F^*}{\F^*}{\I_C}$$
we have 
$$\cdots\rightarrow\Ext^2_{\O_X}(\I_C,\det\F^*)\rightarrow \Ext^2_{\O_X}(\I_C,\F^*)\rightarrow \Ext^2_{\O_X}(\I_C,\I_C)\rightarrow\cdots$$
where $\Ext^2_{\O_X}(\I_C,\det\F^*)=\H^1(X,\I_C\tensor\det\F\tensor\omega_X)=0$ by hypothesis.  Therefore, if $\Ext^2_{\O_X}(\I_C,\I_C)=0$ then $\Ext^2_{\O_X}(\F,\F)=0$.

\underline{Part (2):}

Applying $\Hom_{\O_X}(\cdot,\I_C)$ to the sequence $$\ses{\det\F^*}{\F^*}{\I_C}$$ we have
$$\cdots\rightarrow \Ext^1_{\O_X}(\det\F^*,\I_C)\rightarrow\Ext^2_{\O_X}(\I_C,\I_C)\rightarrow\Ext^2_{\O_X}(\F^*,\I_C) \rightarrow\cdots$$
where $\Ext^1_{\O_X}(\det\F^*,\I_C)=\H^1(X,\det\F\tensor\I_C)=0$ by hypothesis.

Applying $\Hom_{\O_X}(\F^*,\cdot)$ to the sequence $$\ses{\det\F^*}{\F^*}{\I_C}$$ we have
$$\cdots\rightarrow \Ext^2_{\O_X}(\F^*,\F^*)\rightarrow\Ext^2_{\O_X}(\F^*,\I_C)\rightarrow\Ext^3_{\O_X}(\F^*,\det\F^*) \rightarrow\cdots$$
where $\Ext^3_{\O_X}(\F^*,\det\F^*)=\H^3(X,\F^*)=0$ by hypothesis.  Therefore, the vanishing of $\Ext^2_{\O_X}(\F,\F)$ implies $\Ext^2_{\O_X}(\I_C,\I_C)=0$.
\nopagebreak \hfill $\Box$ \par \medskip
\end{proof}

Especially when $X$ is Fano, part (2) of Proposition~\ref{betterhilbvsreflex} is at least consistent with the idea that the Hilbert scheme of curves should fiber over the moduli of reflexive sheaves.  This will be discussed further in Section~\ref{globaltheory}.

If our primary interest is the moduli of reflexive sheaves, then combining Corollary~\ref{F} with part (1) of Proposition~\ref{betterhilbvsreflex} yields the following local result:

\begin{thm}\label{bigthm}
Let $\F$ be a stable reflexive sheaf of rank $2$ on a smooth projective threefold $X$ with $\H^2(X,\F)=0$ and assume that either $\omega_X^*$ is effective or that there exists an $n\in\Z$ such that $\H^0(X,\F\tensor\omega_X^n)\neq0$ and $\H^0(X,\F\tensor\omega_X^{n+1})=0$. 
Suppose further that $\F$ has a section whose zero scheme is a curve $C$, and that
\begin{enumerate}
\item $\H^1(X,\I_C\tensor\det\F\tensor\omega_X)=0$
\item $\Ext^2_{\O_X}(\I_C,\I_C)=0$
\end{enumerate}
Then the (coarse) projective moduli space of semi-stable coherent rank $2$ torsion-free sheaves at the point corresponding to $\F$ is smooth of dimension $\mathfrak{D}(\F)$.
\nopagebreak \hfill $\Box$ \par \medskip
\end{thm}

A particularly nice application of Theorem~\ref{bigthm} is the case $C$ is a rational curve on a Fano threefold.

\begin{cor}\label{fanocor}
Let $X$ be a smooth projective Fano threefold, $\F$ a stable rank $2$ reflexive sheaf with $\det\F$ big and nef.  Suppose that $\F$ has a section whose zero scheme is a rational curve $C$. If $\Ext^1_{\O_X}(\I_C,\O_C)=0$ then the (coarse) projective moduli space of semi-stable coherent rank $2$ torsion-free sheaves is smooth of dimension $$1-\frac{c_1(X)c_2(X)}{6}+\frac{c_1(X)\Delta(\F)}{2}$$ at the point corresponding to $\F$.
\end{cor}

\begin{proof}
By the nef hypothesis, $\H^1(C,\det\F\tensor\O_C)=0$.

From the sequence $$\ses{\det\F^*}{\F^*}{\I_C}$$ we have $\H^2(X,\det\F^*)=0$ by hypothesis and $\H^2(X,\I_C)=0$ because $C$ is rational; hence $\H^2(X,\F^*)=0$.  By Proposition~\ref{localtoglobal}, this gives $\Ext^2_{\O_X}(\F,\O_X)=0$, but $\Ext^2_{\O_X}(\F,\O_X)^*=\H^1(X,\F\tensor\omega_X)$.  From the sequence $$\ses{\omega_X}{\F\tensor\omega_X}{\I_C\tensor\det\F\tensor\omega_X}$$
we have $\H^1(X,\I_C\tensor\det\F\tensor\omega_X)=0$.

From the long exact sequence
$$\cdots\rightarrow\Ext^1_{\O_X}(\I_C,\O_C)\rightarrow \Ext^2_{\O_X}(\I_C,\I_C)\rightarrow \Ext^2_{\O_X}(\I_C,\O_X)\rightarrow \cdots$$
and the fact that $\Ext^2_{\O_X}(\I_C,\O_X)^*=\H^1(X,\I_C\tensor\omega_X)=0$ we have that $\Ext^2_{\O_X}(\I_C,\I_C)=0$.
\nopagebreak \hfill $\Box$ \par \medskip
\end{proof}

\begin{rem}{Remark}{CYcor}
In the case of a canonically trivial threefold, a stable reflexive sheaf cannot have a section whose zero scheme is a rational curve.  However, one can still show that for any reflexive sheaf $\F$ with a section whose zero scheme is an \textit{irreducible} rational curve $C$, if $\H^1(C,N_{C/X})=0$ then $\Ext^2_{\O_X}(\F,\F)=0$.  
\nopagebreak \hfill $\Box$ \par \medskip
\end{rem}

\section{The Global Structure of Moduli of Locally Free Sheaves and Donaldson-Thomas Invariants}\label{globaltheory}

Note that in Theorem~\ref{bigthm} we have the condition that the moduli space of ideal sheaves is smooth at the point $[\I_C]$.  From the point of view of the Serre Correspondence some relationship between the moduli of reflexive sheaves and the moduli of ideal sheaves is not surprising; similarly, we expect a relationship with the Hilbert scheme.  Note, however, that a single reflexive sheaf may have many independent sections whose zero schemes may be different curves within the same component of the Hilbert scheme.  We therefore expect that components of the Hilbert scheme of curves fiber over moduli spaces of reflexive sheaves, with fibers projectivized spaces of sections (since scalar multiples of sections have the same zero scheme).   This phenomenon was already described in \cite{IM2} (see Example~\ref{semicanonical} below).  It should also be noted that it is possible, a priori, that $\F$ could have many independent sections vanishing along the \textbf{same} scheme, though it is proved in \cite[1.3]{hartvb} (the result is attributed to Wever) that this is generally not the case for vector bundles on $\P^3$.  In general, we have

\begin{prop}
Let $X$ be a smooth projective threefold, $\F$ a simple rank $2$ locally free sheaf with $\H^1(X,\det\F^*)=0$.  If $\F$ has a section whose zero scheme is a curve $C$ and if $\H^1(X,\I_C)=0$, then $\h^0(X,\I_C\tensor\F)=1$.
\end{prop}

\begin{proof}
Tensoring the exact sequence
$$\ses{\det\F^*}{\F^*}{\I_C}$$
by $\F$, by simplicity it is enough to show $\H^1(X,\F^*)=0$.  This follows immediately from the hypotheses.
\nopagebreak \hfill $\Box$ \par \medskip
\end{proof}

As evidence supporting these observations (in addition to Examples~\ref{luckwithr=2} and ~\ref{luckwith(2,3)}) we have the following elementary result. 

\begin{prop}\label{hilbdominates}
Let $X$ be a smooth projective threefold, $\F$ a rank $2$ locally free sheaf.   Suppose that $\F$ has a section whose zero scheme is a curve $C$ and suppose that $\H^1(X,\F)=\H^1(X,\F^*)=\H^2(X,\F)=\H^2(X,\F^*)=0$.  Then $\H^0(C,N_{C/X})\twoheadrightarrow\Ext^1_{\O_X}(\F,\F)$ and $\Ext^2_{\O_X}(\F,\F)\hookrightarrow\H^1(C,N_{C/X})$.  In particular, we have
$$\dim\Ext^1_{\O_X}(\F,\F)=\h^0(C,N_{C/X})-\h^0(X,\F)+\h^0(X,\I_C\tensor\F).$$
If, in addition, we have $\H^1(X,\I_C)=0$ then
$$\dim\Ext^1_{\O_X}(\F,\F)=\h^0(C,N_{C/X})-\left(\h^0(X,\F)-1\right).$$
\end{prop}

This gives precisely the naive dimension count of the above discussion.  Compare this with Proposition~\ref{betterhilbvsreflex} and with \cite[6.6]{Vakil} where, in particular, it is shown that if $C\subset X$ is a local complete intersection curve in a smooth threefold (which is the case in Proposition~\ref{hilbdominates}) with $\H^1(X,\O_X)=\H^2(X,\O_X)=0$ then $\H^0(C,N_{C/X})\cong\Ext^1_{\O_X}(\I_C,\I_C)$ and $\H^1(C,N_{C/X})\hookrightarrow\Ext^2_{\O_X}(\I_C,\I_C)$.

\begin{proof}(of Proposition~\ref{hilbdominates})
From the section of $\F$ we have 
$$\ses{\det\F^*}{\F^*}{\I_C}$$
Tensoring with $\F$ and taking global sections gives
$$\cdots\rightarrow\H^i(X,\F^*)\rightarrow\Ext_{\O_X}^i(\F,\F)\rightarrow\H^i(X,\F\tensor\I_C)\rightarrow\cdots$$
Tensoring the standard sequence
$$\ses{\I_C}{\O_X}{\O_C}$$
with $\F$ and taking global sections gives
$$\cdots\rightarrow\H^i(X,\F\tensor\I_C)\rightarrow\H^i(X,\F)\rightarrow\H^i(C,N_{C/X})\rightarrow\cdots$$
The results follow from simple diagram chasing.
\nopagebreak \hfill $\Box$ \par \medskip
\end{proof}

\begin{rem}{Remark}{aCM}
If $\picard{X}=\Z$, then Proposition~\ref{hilbdominates} applies in particular to aCM vector bundles \cite{Beauville} (these are rank $2$ bundles $\F$ such that $\h^i(X,\F(n))=0$ for $i=1,2$ and $n\in\Z$).  See Example~\ref{semicanonical} below.  
\nopagebreak \hfill $\Box$ \par \medskip
\end{rem}

\begin{rem}{Remark}{moreevidence}
Part (2) of Proposition~\ref{betterhilbvsreflex} should also be evidence of this correspondence; i.e. if the Hilbert scheme of curves does, in fact, fiber over the moduli of reflexive sheaves then one may expect that smoothness of the base, together with some extra conditions, would give smoothness of the total space. 
\nopagebreak \hfill $\Box$ \par \medskip
\end{rem}

We do not review Donaldson-Thomas invariants here (nice introductions can, for example, be found in \cite{Katz},\cite{klq},\cite{lq1},\cite{lq2}) but they are a motivation for studying the global relationship between the moduli space of reflexive sheaves and the Hilbert scheme of curves.

For example, fix a subcanonical curve $C$ on a Calabi-Yau threefold.  The expected dimension, and hence the dimension of the virtual fundamental class, of both the component of the Hilbert scheme and the relevant moduli space of torsion-free sheaves is zero.  However, Proposition~\ref{hilbdominates} suggests that the moduli space of torsion-free sheaves may have strictly smaller dimension and so the geometry of that space may be more directly related to invariants coming from virtual fundamental classes.

\begin{rem}{Example}{secants}
Let $C\subset\P^3$ be the twisted cubic, let $X=\operatorname{Bl}_C\P^3$ be the blow up along $C$, and let $L\subset X$ be the proper transform of a secant line to $C$.  Then $\omega_L\tensor\omega_X^*=\O_L$, hence $L$ is subcanonical.  As $\H^1(X,\O(-2H+E))=\H^2(X,\O(-2H+E))=0$ (Kodaira Vanishing) we may associate to $L$ a unique vector bundle $\F$ with $\det\F=\O_X(2H-E)$ (note that $\O_L(2H-E)=\O_L$).

The linear system $|2H-E|$ gives a morphism $f:X\rightarrow\P^2$.  It turns out (\cite{vermeire}) that this is a $\P^1$-bundle; in fact $X=\P(f_*\O_X(H))$.  Given this, it is not hard to see that $\F=f^*f_*\O_X(H)$ and that $\H^1(X,\F)=\H^1(X,\F^*)=\H^2(X,\F)=\H^2(X,\F^*)=0$.  Further, as $L$ is a fiber of $f$, $\h^0(L,N_{L/X})=2$ and $\h^1(L,N_{L/X})=0$.  Finally, the exact sequence
$$\ses{\O_X}{\F}{\I_L(2H-E)}$$
implies that $\h^0(X,\F)=3$.  By Proposition~\ref{hilbdominates}, we have $\Ext^1_{\O_X}(\F,\F)=\Ext^2_{\O_X}(\F,\F)=0$.
\nopagebreak \hfill $\Box$ \par \medskip
\end{rem}

\begin{rem}{Example}{itdoeshappen}
Let $C\subset\P^3$ be a smooth canonical curve of genus $4$.  Every such $C$ lies (degenerately) on a smooth quintic $Q\subset\P^4$.  $C$ is subcanonical since $\omega_C=\O_C(1)$ and $\omega_Q=\O_Q$, and so by the Serre correspondence there is a rank $2$ vector bundle $\F$ on $Q$ with $\det\F=\O_Q(1)$ and with a section whose zero scheme is $C$ (it is easy to see $\F$ is stable, Cf. \cite[3.1]{hartvb}).  The corresponding extension
$$\ses{\O_Q}{\F}{\I_{C/Q}(1)}$$
immediately gives $\h^0(Q,\F)=2$ and $\h^i(Q,\F)=0$ for $i>0$, and dually gives $\h^3(Q,\F^*)=2$ and $\h^i(Q,\F^*)=0$ for $i<3$.  Therefore, Proposition~\ref{hilbdominates} applies and we see that 
$$\dim\Ext^1_{\O_Q}(\F,\F)=\h^0(C,N_{C/Q})-1.$$
\nopagebreak \hfill $\Box$ \par \medskip
\end{rem}

Again, we point to Examples~\ref{luckwithr=2} and ~\ref{luckwith(2,3)} for examples on Fano threefolds.

To construct a map from the component of the Hilbert scheme $\mathcal{H}$ containing a curve $C$ to the component of the moduli space of reflexive sheaves $\mathcal{M}$ containing a corresponding semistable reflexive sheaf $\F$, we interpret the Serre correspondence as a correspondence in $\mathcal{H}\times\mathcal{M}$.  Specifically, consider
\begin{center}
{\begin{minipage}{1.5in}
\diagram
 \Gamma\subset\mathcal{H}\times\mathcal{M} \dto_{\pi_1} \drto^{\pi_2} & \\
 \mathcal{H} & \mathcal{M}
\enddiagram
\end{minipage}}
\end{center} 
where $\Gamma=\{ (C,\F)|\text{ there exists }s\in\Gamma(X,\F)\text{ such that }Z(s)=C\}$.
\begin{prop}\label{serremap}
Notation as above, assume that $C\subset X$ is a connected, subcanonical curve corresponding to a semistable vector bundle $\F$.
Then $\pi_2\circ\pi_1^{-1}(C)$ is a single point $[\F]$; therefore, the induced rational map $\Sigma:\mathcal{H}\dashrightarrow\mathcal{M}$ is defined at $C$.

Suppose additionally that $\H^0(X,\F(-D))=0$ for every non-trivial effective divisor $D$.  Then as a set $\Sigma^{-1}(\F)=\P\Gamma(X,\F)$.
\end{prop}

\begin{proof}
Suppose $\L\tensor\omega_X\tensor\O_C\cong\omega_C$.
Referring to the proof of \cite[1.1]{hartvb}, the main point is that since $\ext^1_{\O_X}(\I_C,\L^*)=1$, there is only one vector bundle $\F$ such that there is a pair $(\F,s)$ with $Z(s)=C$. 

The second part guarantees that every section of $\F$ vanishes along a curve.
\nopagebreak \hfill $\Box$ \par \medskip
\end{proof}

\begin{rem}{Example}{semicanonical}
In \cite[2.4]{IM2} it is shown that if $C$ is an ACM half-canonical curve of genus $15$ on a quartic threefold not contained in a quadric, then the Serre Correspondence gives a morphism from a subscheme of the Hilbert scheme to an open set of the appropriate moduli space of rank $2$ vector bundles.  Further, the fibers are identified with projectivized spaces of sections as in Proposition~\ref{serremap}.
\nopagebreak \hfill $\Box$ \par \medskip
\end{rem}

To illustrate Proposition~\ref{serremap}, we give a simple corollary. 

\begin{cor}\label{birational}
Let $X\subset\P^4$ be a smooth projective hypersurface of degree $d$, $C\subset X$ a connected, nondegenerate, subcanonical curve with $\O_C(d-4)\cong\omega_C$.  Let $\F$ be the corresponding (necessarily stable) vector bundle with $\det\F=\O_X(1)$.  Then the reduced component of the Hilbert scheme containing $[C]$ is birationally equivalent to the reduced component of the moduli space of rank $2$ vector bundles containing $\F$.
\nopagebreak \hfill $\Box$ \par \medskip
\end{cor}

\begin{rem}{Example}{someexamplesattheend}
Some examples of Corollary~\ref{birational} include:
\begin{enumerate}
\item $C$ is an elliptic normal of degree $5$ curve lying on a quartic.
\item $C$ is a canonical curve of genus $5$ lying on a quintic.
\item $C$ is a nondegenerate half-canonical curve lying on a sextic (e.g. of genus $15$ \cite{IM2}).
\end{enumerate}
\nopagebreak \hfill $\Box$ \par \medskip
\end{rem}

{\bf Acknowledgments:}
I would like to thank Robin Hartshorne for an enlightening exchange; specifically, the proof of Proposition~\ref{killext3fano} is essentially due to him and is much better than the somewhat awkward proof I had originally.  I would also like to thank Chad Schoen for some general discussions and Aaron Bertram for suggesting that I consider virtual fundamental classes.

\end{document}